\documentclass[12pt]{article}
\usepackage{amsmath}
\usepackage{amssymb}
\usepackage{amscd}
\usepackage{latexsym}
\usepackage{amsthm}
\usepackage{setspace}
\usepackage{epsfig}
\usepackage{psfrag}
\usepackage{amsfonts}
\usepackage{enumerate}
\usepackage{xypic}

\newtheorem{theorem}{Theorem}[section]
\newtheorem*{theorem*}{Theorem}
\newtheorem{lemma}[theorem]{Lemma}
\newtheorem{proposition}[theorem]{Proposition}
\newtheorem{corollary}[theorem]{Corollary}

\newtheorem{conjecture}[theorem]{Conjecture}

{
\theoremstyle{plain}

\newtheorem{example}[theorem]{Example}

\newtheorem{remark}[theorem]{Remark}
}

\newcommand{\Hom}{\operatorname{Hom}}
\newcommand{\Char}{\operatorname{Char}}
\newcommand{\Stab}{\operatorname{Stab}}
\newcommand{\Vol}{\operatorname{Vol}}

\renewcommand{\Im}{\operatorname{Im}}

\newcommand{\Spec}{\operatorname{Spec}}

\renewcommand{\dim}{\operatorname{dim}}
\newcommand{\Ann}{\operatorname{Ann}}
\newcommand{\Sym}{\operatorname{Sym}}

\newcommand{\Od}{\displaystyle\bigoplus}

\newcommand{\C}{{\mathbb{C}}}
\newcommand{\Z}{{\mathbb{Z}}}
\newcommand{\Q}{{\mathbb{Q}}}

\newcommand{\R}{{\mathbb{R}}}

\newcommand{\otn}{\{1,\ldots,n\}}
\newcommand{\half}{\frac{1}{2}}
\newcommand{\bd}{\partial}
\newcommand{\la}{\lambda}

\newcommand{\subs}{\subseteq}

\newcommand{\hookto}{{\hookrightarrow}}

\newcommand{\D}{\Delta}
\renewcommand{\i}{\iota}
\renewcommand{\a}{\alpha}

\renewcommand{\t}{\theta}

\newcommand{\eps}{\varepsilon}

\renewcommand{\cot}{T^*V}

\newcommand{\Rn}{\R^n}

\newcommand{\minv}{\mu^{\! -1}(\la)}

\renewcommand{\mod}{{\!/\!\!/\!}}

\newcommand{\Tn}{T^n}
\newcommand{\Td}{T^d}

\newcommand{\tn}{\mathfrak{t}^n}

\newcommand{\td}{\mathfrak{t}^d}

\newcommand{\tnd}{(\tn)^*}

\newcommand{\tr}{\t_{\R}}
\newcommand{\tnr}{\tn_{\R}}
\newcommand{\tdr}{\td_{\R}}
\newcommand{\tndr}{(\tn_{\R})^*}
\newcommand{\tddr}{(\td_{\R})^*}

\renewcommand{\t}{\mathfrak{g}}
\newcommand{\tdu}{\t^*}
\newcommand{\stab}{\mathfrak{stab}}

\newcommand{\vst}{V^{st}}
\renewcommand{\vss}{V^{ss}}
\newcommand{\vlf}{V^{\ell\! f}}
\newcommand{\evp}{\operatorname{ev}(p)}

\newcommand{\Gic}{G_i^\chi}
\newcommand{\Fic}{F_i^\chi}

\newcommand{\PrA}{\Pr_{\! A}}
\renewcommand{\Pr}{P^{\chi}}
\newcommand{\dc}{\D^{\chi}}
\newcommand{\dca}{\dc_A}
\newcommand{\Dr}{\dc}
\newcommand{\DrA}{\dca}

\newcommand{\rt}{\tilde \chi}
\newcommand{\odra}{\mathbf{1}_{\DrA}}

\newcommand{\WA}{W_A}
\newcommand{\WAbd}{W_A^{bd}}
\newcommand{\Fbd}{\mathcal{F}^{bd}}
\newcommand{\F}{\mathcal{F}}

\renewcommand{\(}{\left(}
\renewcommand{\)}{\right)}

\renewenvironment{proof}{\noindent {\bf Proof:}}{\qed \par}
\newenvironment{proofint}{\noindent {\bf Proof of \ref{int}:}}{\qed \par}

\begin{document}
\setstretch{1.05}
\baselineskip=\baselineskip

\noindent {\LARGE \bf All the GIT quotients at once}
\bigskip\\
{\bf Nicholas Proudfoot}\footnote{Partially supported
by a National Science Foundation Postdoctoral Research Fellowship.}\\
Department of Mathematics, University of Texas,
Austin, TX 78712
\bigskip
{\small
\begin{quote}
\noindent {\em Abstract.}
Let $G$ be an algebraic torus acting on a smooth variety $V$.
We study the relationship between the various GIT quotients of $V$
and the symplectic quotient of the cotangent bundle of $V$.
\end{quote}
}
\bigskip

\noindent
Let $G$ be a reductive algebraic group acting on a smooth variety $V$.
The cotangent bundle $T^*V$ admits a canonical algebraic
symplectic structure, and the induced action of $G$ on $T^*V$ is hamiltonian, that is,
it admits a natural moment map $\mu:T^*V\to\mathfrak g^*$ (see Equation \eqref{mm}
for an explicit formula).  Over the past ten years, a guiding principle has emerged that
says that if $X$ is an interesting variety which may be naturally presented as a GIT
(geometric invariant theory) quotient of $V$ by $G$, then the symplectic quotient
$\mu^{-1}(\la)\,\mod\, G$ of $T^*V$ by $G$ is also interesting.
This mantra has been particularly fruitful on the level of cohomology, as we describe below.
Over the complex numbers, a GIT quotient may
often be interpreted as a K\"ahler quotient by the compact form of $G$, and an algebraic
quotient as a hyperk\"ahler quotient.  For this reason, the symplectic quotient may be
loosely thought of as a quaternionic or hyperk\"ahler analogue of $X$.  
Let us review a few examples of this construction.

\vspace{.2cm}
\noindent{\bf Hypertoric varieties.}
These examples comprise the case where $G$ is abelian and $V$ is a linear
representation of $G$.  
The geometry of toric varieties is deeply related to the combinatorics of polytopes;
for example, Stanley \cite{St} used the hard Lefschetz theorem for toric varieties
to prove certain inequalities for the $h$-numbers of a simplicial polytope.
The hyperk\"ahler analogues of toric varieties, known as hypertoric varieties,
interact in a similar way with the combinatorics of rational hyperplane arrangements.
Introduced by Bielawski and Dancer \cite{BD}, hypertoric varieties were used by
Hausel and Sturmfels \cite{HS} to give a geometric interpretation of virtually every known property
of the $h$-numbers of a rationally representable matroid.  Webster and the author \cite{PW}
extended this line of research by studying the intersection cohomology groups of singular
hypertoric varieties.

\vspace{.2cm}
\noindent{\bf Quiver varieties.}
A quiver is a directed graph, and a representation of a quiver is a vector
space for each node along with a linear map for each edge.
For any quiver,
Nakajima \cite{N1,N2,N3} defined a quiver variety
to be the quaternionic analogue of the moduli space of framed representations.
Examples include the Hilbert scheme of $n$ points in the plane and the moduli space
of instantons on an ALE space.  He has shown that the cohomology and K-theory
groups of quiver varieties carry actions of Ka\v{c}-Moody algebras and their associated
Hecke algebras, and has exploited this fact to define canonical bases for highest weight
representations.  Crawley-Boevey and Van den Bergh \cite{CBVdB} and Hausel \cite{Ha}
have used Betti numbers of quiver varieties to prove a long standing conjecture of Ka\v{c}.
 
\vspace{.2cm}
\noindent{\bf Hyperpolygon spaces.}
Given an ordered $n$-tuple of positive real numbers, the associated polygon space
is the moduli space of $n$-sided polygons in $\R^3$ with edges of the prescribed lengths,
up to rotation.  Such a space may be interpreted as a moduli space of stable configurations
of points on a projective line, or, via the Gelfand-MacPherson correspondence \cite{HK},
as a GIT quotient of the grassmannian $G(2,n)$ by the natural action of the torus $T^{n-1}$.
The quaternionic analogues of these spaces were introduced by Konno \cite{K2}, and dubbed hyperpolygon spaces in \cite{HP}.
In \cite{HP}, Harada and the author show that certain projective
subvarieties of a hyperpolygon space have interpretations in terms of spacial
polygons, which suggests the general problem of searching for moduli-theoretic interpretations
of hyperk\"ahler analogues of moduli spaces.
While hyperpolygon spaces are in fact special cases of quiver varieties, they will be of special
interest to us in this paper because they may be constructed using an abelian group.

\vspace{.2cm}
To define a GIT quotient of a variety $V$ by a group $G$ one needs more data than just an action;
to be precise, we need a $G$-equivariant ample line bundle on $V$.  If we define two
such line bundles to be equivalent whenever they lead to the same GIT quotient, there
will in general be finitely many distinct equivalent classes of equivariant ample line bundles.
In the toric case, these classes correspond to triangulations of an oriented matroid \cite{Sa}.
In the case of polygon spaces, the choice of line bundle corresponds to the choice of edge lengths.
The various GIT quotients of $V$ by $G$ will always be birational, but will generally be topologically
distinct.  

The symplectic quotient of $T^*V$ by $G$ requires two choices, namely an equivariant
line bundle as well as an element $\la\in(\mathfrak{g}^*)^G$ at which to reduce.  If $G$
is abelian and $\la$ is chosen generically, however, then $G$ will act locally freely on $\mu^{-1}(\la)$,
and the symplectic quotient 
$$M_\la := \mu^{-1}(\la)\,\mod\, G = \mu^{-1}(\la)/G$$ will not depend on the choice
of line bundle.  In fact, the topological type of the symplectic quotient over the complex numbers
does not depend on the choice of generic $\la$, either!  Intuitively, this can be seen from the
fact that the set of generic parameters $\la$ is connected, but the noncompactness of the quotients
makes this argument technically difficult.  In this paper we take a different approach, proving
the following theorem (Corollary \ref{all}).

\begin{theorem*}  If $G$ is abelian and $\la\in\mathfrak g^*$ is a regular value for $\mu$,
then $M_\la$ is isomorphic to the total space of an affine bundle
over the nonseparated prevariety obtained by gluing together the various smooth GIT quotients
of $V$ along the open sets on which they agree.
\end{theorem*}

\noindent This surprising result, taken over the complex numbers, implies 
that the topology of any one symplectic quotient of $T^*V$
is intimately related to the topology of {\em all} of the different GIT quotients of $V$.  It is for this
reason that we use the phrase `all the GIT quotients at once'.  We note that the theorem
in fact holds over arbitrary fields, and was used in \cite{PW} to count points on hypertoric varieties
over finite fields.

Section \ref{git} is devoted to giving a careful definition of the various objects referred to in the above
theorem.  The proof itself is remarkably simple, following Crawley-Boevey's work on quiver
varieties in \cite{CB}.  In Section \ref{conj} we consider the natural map from the cohomology ring
of $M_\la$ to the direct sum of the cohomology rings of the GIT quotients, whose existence
follows immediately from the theorem.  Konno \cite[7.6]{K2} 
proves that this map is an injection in the case
of hyperpolygon spaces\footnote{Konno does not phrase his theorem in these terms,
as he does not have Corollary \ref{all} at his disposal.  But it is easy to translate his result
into the one that we attribute to him.}, and we prove the analogous theorem for hypertoric varieties 
(Theorem \ref{int}).  We conjecture that such a result 
will hold in greater generality (Conjecture \ref{conj}).  A consequence of this conjecture would
be that the cohomology ring of $M_\la$ is level, meaning that it satisfies an analogue of
Poincar\'e duality for noncompact spaces (see Remark \ref{level}).
In the case of hypertoric varieties, this is a well known and nontrivial fact \cite[\S 7]{HS},
but it is by no means the case for every smooth manifold.

\begin{section}{An affine bundle}\label{git}
Let $V$ be a smooth algebraic variety over an arbitrary field $k$.
We will assume that $V$ is projective over affine, which means that the natural
map $V \to \Spec \mathcal{O}_V$ is projective.
Let $G$ be an algebraic torus over $k$ acting effectively on $V$, and let $L$ be a $G$-equivariant
ample line bundle on $V$.  A point $p\in V$ is called {\em $L$-semistable}
if there exists a $G$-invariant section of a positive power of $L$ that does not vanish at $p$.
The set of $L$-semistable points of $V$ will be denoted $V^{ss}(L)$.
An $L$-semistable point $p$ is called {\em $L$-stable} if $G$ acts locally freely at $p$
and its orbit is closed in $V^{ss}(L)$.  The set of $L$-stable points of $V$
will be denoted $\vst(L)$.  If every $L$-semistable point is $L$-stable,
then we will call $L$ {\em nice}.

We will consider two equivariant ample line bundles to be equivalent if they induce
the same stable and semistable sets.  Let $\{L_i\mid i\in I\}$ be a complete set of representatives
of equivalence classes of nice line bundles with nonempty stable sets.
Let $\vlf$ denote the set of points of $V$ at which $G$ acts locally freely.  
By definition, $\vst(L)$ is contained in $\vlf$ for any $L$.
The following lemma is a converse to this fact.

\begin{lemma}\label{lf}
Suppose that there exists at least one ample equivariant line bundle on $V$.
If $G$ acts locally freely at $p$, then $p$ is $L$-stable for some nice $L$,
thus we have $\vlf = \displaystyle\bigcup_{i\in I}\vst(L_i)$.
\end{lemma}

\begin{proof}
For any point $p\in V$,
and any equivariant line bundle $L$, 
choosing an identification of $L_p$ with $k$
gives us a natural element $\evp\in\Gamma(L)^\vee$, the dual of the vector
space of sections of $L$.  
The equivariant structure on $L$ gives a decomposition
$$\Gamma(L)^\vee = \bigoplus_\chi\Gamma(L)^\vee_\chi,$$
where $\chi$ ranges over the lattice $\Char(G) = \Hom(G,\mathbb{G}_m)$.
Let $\evp_\chi$ be the component of $\evp$ corresponding to the
character $\chi$.
The {\em state polyhedron}
$\D_p(L) \subs\Char(G)_\Q$ is defined to be the convex hull inside of $\Char(G)_\Q$
of the set $\{\chi\mid \evp_\chi\neq 0\}$.
Dolgachev and Hu \cite[1.1.5]{DH} show that $p$ is $L$-semistable
if and only if the trivial character is contained in
$\D_p(L)$, and $L$-stable if and only if it is contained in the
interior of $\D_p(L)$.
The affine span of $\D_p(L)$ is equal to a shift of the set
of rational characters that vanish on the stabilizer of $p$, where the
shift is given by the character with which the stabilizer of $p$ acts on $L_p$.
In particular, the interior of $\D_p(L)$ is nonempty if and only if
$G$ acts locally freely at $p$.  

Let $L_\chi$ be the trivial bundle
on $V$ with equivariant structure given by the character $\chi$.
Then $$\D_p(L^{\otimes m}\otimes L_\chi)
= m\cdot \D_p(L) - \chi,$$
thus by taking a high enough tensor power of $L$ and twisting it 
by an appropriate character, we may find a new equivariant ample 
line bundle whose state polytope
is an arbitrary dilation and translation of that of $L$.  
In particular, if $p\in\vlf$, we may find an $L$ with respect to which $p$ is stable.
It remains to show that this line bundle can be chosen to be nice.
For all $q\notin\vlf$, $\D_q(L)$ is contained in a proper affine
subspace of $\Char(G)_\Q$.  Let $\chi$ be a character which
is nontrivial on the stabilizer of $q$ for {\em every} $q\in V\smallsetminus\vlf$.
By again replacing $L$ with a large tensor power and
twisting by $\chi$, we may ensure that none of these subspaces
contains the origin.  If $\chi$ is chosen to be small with respect
to the size of the tensor power, then this operation will not
break the $L$-stability of $p$.
\end{proof}

\vspace{\baselineskip}
For any ample equivariant line bundle $L$, the GIT (geometric invariant
theory) quotient
$$V\mod_L G := \vss(L)/ G$$
is defined to be the categorical quotient of $\vss(L)$ by $G$, 
in which two points are identified if the closures of their $G$-orbits intersect.
The fact that $V$ is projective over $\Spec \mathcal{O}_V$ implies that $V\mod_L G$
is projective over $\Spec\mathcal{O}_V^G$.
If $L$ is nice, then $V\mod_L G$ is simply the geometric quotient of $\vst(L)$ by $G$.
For all $i\in I$, let $X_i = V\mod_{L_i} G$ be the corresponding GIT quotient.

\begin{example}\label{polygons}\em
Let $V = G(2,n)$ be the grassmannian of $2$-planes in $\C^n$, and let $G = T^n/\C^\times_{diag}$
be the $(n-1)$-torus acting on $V$.  The GIT quotients $\{X_i\}$ can also be realized as GIT
quotients of an $n$-fold product of projective lines by the diagonal action of $PSL(2,\C)$.  These
quotients have been studied extensively, for example in {\cite{MFK,Th}}.  
In the symplectic geometry
literature these spaces are known as {\em polygon spaces}, as they parameterize $n$-sided
polygons in $\R^3$ with fixed edge lengths, up to rotation {\cite{HK}}.
\end{example}

The action of $G$ on $V$ induces a symplectic action on the cotangent bundle
$\cot$, with moment map $\mu:\cot\to\tdu$ given by the equation
\begin{equation}\label{mm}
\mu(p,\a)(x) = \a(\hat x_p),
\end{equation}
where $\a$ is a cotangent vector to $V$ at $p$, and $\hat x_p$ is the tangent 
vector at $p$ induced by the infinitesimal action of $x\in\t$.
For $\la\in\tdu$, let $\phi_\la:\minv\to V$ be the canonical projection to $V$.
Given an equivariant ample line bundle $L$ on $V$, we define the space
$$M_{\la,L} := \minv\,\mod_{\phi_\la^*L} G$$ to be the GIT quotient of $\minv$ by $G$.
If $\lambda$ is a regular value of $\mu$, then $G$ acts locally freely on $\minv$,
and the quotient is geometric; in particular, it is independent of $L$.
Thus when $\la$ is a regular value we will drop $L$ from the notation.

\begin{proposition}\label{image}
If $\la\in\tdu$ is a regular value of $\mu$, then $\Im\phi_\la = \vlf$,
and the fibers of $\phi_\la$ are affine spaces of dimension $\dim V - \dim G$.
\end{proposition}

\begin{proof}
For any point $p\in V$, we have
an exact sequence of $k$-vector spaces\footnote{The analogous exact sequence
in the context of representations of quivers
first appeared in \cite{CB}, and was used to count points on quiver varieties
over finite fields in \cite{CBVdB}.}
\begin{equation}\label{cb}
0\to\{\alpha\mid\mu(p,\alpha)=0\}\to T_p^*V\overset{\mu(p,-)}{\longrightarrow}
\tdu\to\stab(p)^*\to 0,
\end{equation}
where $\stab(p)^* = \tdu/\stab(p)^\perp$ is the Lie coalgebra of the
stabilizer of $p$ in $G$.
By exactness of \eqref{cb} at $\tdu$, $p$ is in the image of $\phi_\la$ if and only if
$\lambda\cdot\stab(p)=0$.  If $p\in\vlf$, then $\stab(p) = 0$, thus $p\in\Im\phi_\la$.
Conversely, suppose that $p\in\Im\phi_\la$.  Since $\la$ is a regular value of $\mu$,
$G$ acts locally freely on $\minv$, and therefore $\Stab(p)\subs G$ acts locally freely
on $\phi_\la^{-1}(p)$.  But $\phi_\la^{-1}(p)$ is a torsor for the vector space 
$\{\alpha\mid\mu(p,\alpha)=0\}$, and any torus action on an affine space has a fixed point.
Hence $\Stab(p)$ must be finite, and therefore $p\in\vlf$.
Finally, we see that if $p\in\vlf=\Im\phi_\la$, then $\dim\phi_\la^{-1}(p)
= \dim\{\alpha\mid\mu(p,\alpha)=0\}=\dim T^*_pV-\dim\tdu$.
\end{proof}

\vspace{.2cm}
\begin{remark}\em
The space $M_\la$ is sometimes known as the {\em twisted cotangent bundle} to the stack $V/G$.
In this language, Proposition \ref{image} says that the $M_\la$ has support $X^{\ell f}\subs V/G$,
and that over its support it has constant rank.
\end{remark}

Consider the nonseparated prevariety
$$X^{\ell\! f} = \vlf\!/G = \bigcup_{i\in I}\vst(L_i)/G = \bigcup_{i\in I}X_i,$$
where two GIT quotients $X_i$ and $X_j$ are glued together along the open set of
points which are simultaneously stable for both $L_i$ and $L_j$.
Proposition \ref{image} has the following immediate corollary.

\begin{corollary}\label{all}
For any regular value $\la$ of $\mu$,
$M_\la = \minv/G$ is isomorphic to the total space of an affine bundle over $X^{\ell\! f}$,
modeled on the cotangent bundle.
\end{corollary}

\end{section}

\begin{section}{A cohomological conjecture}\label{cohom}
The purpose of this section will be to consider the cohomological implications of
Corollary \ref{all}, taken over the complex numbers,
focusing on the case of toric and hypertoric varieties.
Specifically, there is a natural map on cohomology
\begin{equation}\label{psi}
\Psi:H^*(M_\la;\R)\cong H^*(X^{\ell\! f};\R)\to\Od_{i\in I}H^*(X_i;\R)
\end{equation}
given by the inclusions of each $X_i$ into $X^{\ell\! f}$.

\begin{conjecture}\label{conj}
If $\mathcal{O}_V^G\cong\C$, then $\Psi$ is injective.
\end{conjecture}

\begin{remark}\em
The hypothesis $\mathcal{O}_V^G\cong\C$ says precisely that the GIT quotients $\{X_i\}$
are projective, and without this assumption Conjecture \ref{conj} fails.  Specifically, let $G=\C^\times$
act on $V=\C^2$ with eigenvalues $\pm 1$.  There are two GIT quotients, namely
$X_1 = \big(V\smallsetminus \C\times\{0\}\big)/G\cong \C$ and
$X_2 = \big(V\smallsetminus \{0\}\times\C\big)/G\cong \C$, neither of which has any nontrivial
cohomology.  On the other hand, $X^{\ell\! f} = V\smallsetminus\{0\}/G$ is an affine line
with a double point, which is weakly homotopy equivalent to the 2-sphere.  So $\Psi$
fails to be injective in degree 2.
\end{remark}

\begin{remark}\em
Conjecture \ref{conj} can be immediately deduced in the case of hyperpolygon spaces
(see Example \ref{polygons}) from \cite[7.6]{K2}.
\end{remark}

\begin{remark}\label{level}\em
If Conjecture \ref{conj} were true, it would imply that $H^*(M_\la;\R)$ is {\em level},
which means that every nonzero class divides a nonzero class of top degree.
(One may think of this property as a generalization of Poincar\'e duality to the situation
where the top degree cohomology need not be one dimensional.)  Indeed, given
a nonzero class $\alpha\in H^*(M_\la;\R)$, Conjecture \ref{conj} would tell us that $\alpha$
restricts to a nonzero class $\alpha_i\in H^*(X_i;\R)$ for some $i$.  Then, by Poincar\'e
duality for $X_i$, there exists a class $\beta_i\in H^*(X_i;\R)$ such that $\alpha_i\cdot\beta_i$
is a nonzero class on $X_i$ of top degree.  Now it suffices to show that $\beta_i$ is the 
restriction of a class in $H^*(M_\la;\R)\cong H^*(X^{\ell\! f};\R)$.  But this follows from Kirwan
surjectivity \cite[5.4]{Ki}, which tells us that the map $H_G^*(V;\R)\to H^*(X_i;\R)$ is surjective.
This map factors through $H^*_G(\vlf;\R)\cong H^*(X^{\ell\! f};\R)$, which implies that the map
from $H^*(X^{\ell\! f};\R)$ to $H^*(X_i;\R)$ is surjective as well.

We note that Conjecture \ref{conj} is in fact equivalent to the statement that
$H^*(M_\la;\R)$ is level {\em and} the fundamental cycles of the GIT quotients
$\{X_i\}$ generate the top homology of $X^{\ell\! f}$.
\end{remark}

The rest of Section \ref{cohom} will be devoted to understanding and proving Conjecture
\ref{conj} in the toric case (Theorem \ref{int}).  Let
$V = \C^n$, $\Tn$ the coordinate torus acting on $V$, and $G$ is a 
codimension $d$ subtorus of $\Tn$.
Let $\chi:\Tn\to \C^*$ be a multiplicative character defined by the equation
$$\chi(t) = t_1^{\chi_1}\ldots t_n^{\chi_n},$$
and let $L_\chi$ be the $G$-equivariant line bundle on $V$ obtained from twisting
the trivial bundle by the restriction of $\chi$ to $G$.
The GIT quotients $V\mod_{L_\chi}G$ are called {\em toric varieties}\footnote{An introduction
to toric varieties from the GIT perspective can be found in \cite{Pr}.},
and the symplectic quotients $M_{\la,L_\chi}$ are called {\em hypertoric varieties}.
Both are intricately related to various combinatorial data associated to the subtorus $G\subs T^n$
and the character $\chi$, as we describe below.

Let $\Td = \Tn/G$ with Lie algebra $\td=\tn/\t$.  This algebra is equipped with an integer lattice
(the kernel of the exponential map), and therefore with a canonical real form $\tdr\subs\td$.
Let $a_1,\ldots,a_n\in\tdr$ be the projections of the standard
basis vectors in $\tn$.  
For all $i\in\otn$, we define half spaces
\begin{equation}\label{fg}
\Fic := \{x\in\tddr\mid x\cdot a_i + \chi_i \geq 0\}\hspace{.7cm}
\text{and}\hspace{.7cm}
\Gic := \{x\in\tddr\mid x\cdot a_i + \chi_i \leq 0\}.
\end{equation}
The geometry of the toric variety $V\mod_{L_\chi}G$ is completely controlled by the polyhedron
\begin{equation}\label{dchi}
\D^\chi := \bigcap_{i=1}^n \Fic
\end{equation}
(see \cite{Pr} and references therein).
The line bundle $L_\chi$ is nice if and only if $\D^\chi$ is a simple polyhedron of dimension $d$,
which means that exactly $d$ facets meet at each vertex.
In this case, we have the following well-known theorem (see for example \cite[2.11]{HS}).

\begin{theorem}\label{tco}
If $L_\chi$ is nice, then the $\Td$-equivariant cohomology ring of $V\mod_{L_\chi}G$ is isomorphic to the Stanley-Reisner ring of the normal fan to $\Dr$.
\end{theorem}

In particular, this implies that the Betti numbers of toric varieties are given
by combinatorial invariants of polytopes \cite{St}.
Theorem \ref{tco} has the following analogue for hypertoric varieties, which appeared
in this form in \cite{HS}, and in a different but equivalent form in \cite{K1}.
\begin{theorem}\label{htco}
If $\la$ is generic, then the $\Td$-equivariant cohomology ring of $M_\la$ is isomorphic to the
Stanley-Reisner ring of the matroid associated to the vector collection $\{a_1,\ldots,a_n\}$.
\end{theorem}

A consequence of this fact is that the Betti numbers of hypertoric varieties are combinatorial invariants
of matroids \cite[6.6]{HS}.  It is also possible to obtain this result by counting
points on $M_\la$ over finite fields, as in \cite{PW} or \cite{Ha}.

\begin{remark}\label{omat}\em
The set of equivalence
classes of nice line bundles on $V$ corresponds to the set of possible combinatorial
types of $\Dr$, which in turn is indexed by triangulations of the {\em oriented matroid}
determined by the vectors $\{a_1,\ldots,a_n\}$ \cite{Sa}.
\end{remark}

\begin{theorem}\label{int}
If $\mathcal{O}_V^G\cong \C$, then the map $\Psi$ of Equation \eqref{psi} is injective.
\end{theorem}

\begin{remark}\em
Using Theorems \ref{tco} and \ref{htco}, it is relatively easy to prove injectivity
of the natural lift of $\Psi$ to a map in $\Td$-equivariant cohomology.  Theorem \ref{int}, however,
does not follow formally from this fact.  For our proof it is necessary to use different descriptions
of the cohomology rings of toric and hypertoric varieties (see Theorems \ref{toric} and \ref{hypertoric}).
\end{remark}

\begin{remark}\em
Since we know that the cohomology ring of a hypertoric variety is level \cite[\S 7]{HS},
Remark \ref{level} tells us that Theorem \ref{int} is equivalent to the statement that
the fundamental cycles of the projective toric varieties $\{X_i\}$ 
generate $H^{2d}(X^{\ell\! f})$.  This is a little bit surprising, as $X^{\ell\! f}$ will often
contain proper but nonprojective toric varieties as open subsets.  Theorem \ref{int}
asserts that the fundamental cycle of any such subvariety can be expressed
as a linear combination of the fundamental cycles of the $\{X_i\}$.  Our proof,
though not in this language, will roughly follow this line of attack.
\end{remark}

\begin{proofint}
We begin by noting that the definitions given in Equations \eqref{fg} and \eqref{dchi} make sense
for any vector $\chi = (\chi_1,\ldots,\chi_n)\in\Rn\cong\tndr$, regardless of whether the coordinates of $\chi$ are integers.
For any subset $A\subs\otn$, let
$$\dca\,\, :=\,\, \bigcap_{i\in A}\Fic\,\,\cap\,\,\bigcap_{j\in A^c}G_j^\chi,$$
and note that $\dc_{\emptyset} = \dc$.  
In general, one should imagine $\dca$ as the polytope
obtained from $\dc$ by ``flipping'' it over the hyperplane $\Fic\cap\Gic$
for each $i\in A$.
We will call $\chi$ {\em simple} if for every $A\subs\otn$, $\dca$ is either empty or simple of
dimension $d$.  In particular, if $\chi$ is an integer vector and $\chi$ is simple,
then $L_\chi$ is nice.
If $\DrA$ is nonempty, then it is bounded if and only if the vectors
$\{\eps_1(A)a_1,\ldots,\eps_n(A)a_n\}$ span $\tdr$ over the non-negative real numbers,
where $\eps_i(A) = (-1)^{|A\cap\{i\}|}$.  We call such an $A$ {\em admissible}.
Note that the assumption $\mathcal{O}_V^G = \C$ is equivalent to the assumption
that the empty set is admissible.

The volume function $\Vol\dca$ is locally polynomial in $\chi$.
More precisely, for every simple $\chi\in\tndr$ and every admissible $A\subs\otn$, 
there exists a polynomial $\PrA\in\Sym^d\tnr$ such that for every simple $\eta\in\tndr$ sufficiently
close to $\chi$,
we have $$\Vol\D^{\eta}_A=\PrA(\eta).$$  We will refer to $\PrA$ as the
{\em volume polynomial} of $\DrA$.  The fact that the volume of a polytope
is translation invariant tells us that $\PrA$ lies in the image of the inclusion
$\i:\Sym^d\t\hookto\Sym^d\tn$.  The cohomology rings of toric and hypertoric varieties
may be described in terms of these volume polynomials.

\begin{theorem}\label{toric}{\em\cite{GS,KP}}
If $\chi\in\tndr$ is simple with integer coordinates, then
$$H^*(V\mod_{L_\chi}G;\R) \cong \Sym\tr/\!\Ann(\Pr),$$
where $\Ann(\Pr) = 
\big\{\bd\in\Sym\tdr \mid \bd\cdot(\i^{-1}\Pr) = 0\big\}.$
\end{theorem}

\begin{theorem}\label{hypertoric}{\em\cite[7.1]{HS}}
If $\chi\in\tndr$ is simple with integer coordinates, then for any $\la\in\td$,
$$H^*(M_{\la,L_\chi};\R) \cong \Sym\tr/\!\Ann\{\PrA\mid A \text{ admissible}\}.$$
\end{theorem}

If we fix $\la$ to be a regular value of $\mu$, then we have already observed
that the left-hand side of the isomorphism of Theorem \ref{hypertoric} does not depend
on $\chi$.  From this it follows that the right-hand side does not depend on $\chi$ either;
in other words, the linear span
$$U = \R\big\{\PrA \mid A\text{ admissible}\big\}$$
is independent of $\chi$.  
Since the empty set is admissible, $\Pr$ is contained in $U$ for all simple integral $\chi$,
and this inclusion induces the canonical map from $H^*(M_\la;\R)$ to $H^*(V\mod_{L_\chi}G;\R)$.
The kernel of $\Phi$ is therefore equal to the image in $H^*(M_\la)$ of the set of polynomials
in $\Sym\tdr$ that annihilate $\Pr$ for every simple integral $\chi$.
To prove Theorem \ref{int}, we need to show that every such polynomial annihilates $U$;
in other words, we must show that $U$ is contained in (and therefore equal to) the the linear span
$\R\{\Pr\mid \chi\text{ simple and integral}\}.$

Let $\F$ be the infinite dimensional vector space consisting of all real-valued
functions on $\tddr$, and let $\Fbd$ be the subspace consisting of functions
with bounded support.  For all subsets $A\subs\otn$, let
$$\WA = \Q\big\{\odra\mid \chi\text{ simple and integral}\big\}$$
be the subspace of $\F$ consisting of finite linear combinations
of characteristic functions of polyhedra $\DrA$ for simple and integral $\chi$, and let
$$\WAbd = \WA\cap\Fbd.$$
Note that $\WAbd=\WA$ if and only if $A$ is admissible.

\begin{lemma}\label{key}
For all $A,A'\subs\otn$, $\WAbd = W_{A'}^{bd}$.
\end{lemma}

\begin{proof}
We may immediately reduce to the case where $A' = A\cup\{j\}$.
Fix a simple $\chi\in\tndr$.
Let $\rt\in\tnd$ be another simple element obtained from $\chi$ by putting
$\rt_i=\chi_i$ for all $i\neq j$, and $\rt_j = N$ for some integer $N\gg 0$.
Then $\D_A^\chi\subs\D_A^{\rt}$, and
\begin{eqnarray*}
\D_{A}^{\rt} \smallsetminus \D_{A}^{\chi} &=&
\bigcap_{i\in A}\Fic\,\,\cap\bigcap_{j\in (A')^c} G_j^\chi\,\,\cap\,\,\big(G_k^{\chi}\smallsetminus G_k^{\tilde\chi}\big)\\
&=& \D_{A'}^{\chi} \cap F_k^{\tilde\chi}.
\end{eqnarray*}
Suppose that $f\in\Fbd$ can be written as a linear combination
of functions of the form $\mathbf{1}_{\D_{A'}^{\chi}}$.  Choosing $N$ large enough
that the support of $f$ is contained in the half space 
$F_i^{\tilde\chi} = \{x\mid x\cdot a_i + N\geq 0\}$,
the above computation shows that $f$ can be written
as a linear combination of functions of the form $\odra$,
hence $W_{A'}^{bd}\subs W_{A}^{bd}$.  The reverse inclusion is obtained
by an identical argument.
\end{proof}

By Lemma \ref{key}, we may write
\begin{equation}\label{express}
\odra = \sum_{j=1}^m\a_j \mathbf{1}_{\D^{\eta^j}}
\end{equation}
for any simple $\chi$ and admissible $A$, where $\a_j\in\Z$ and
$\eta^j$ is a simple element of $\tnd$ for all $j\leq m$.
Taking volumes of both sides of the equation, we have
\begin{equation}\label{polys}
\PrA(\chi) = \sum_{j=1}^m\a_j P^{\eta^j}\big(\eta^j\big).
\end{equation}
Furthermore, we observe from the proof of Lemma \ref{key}
that for all $j\leq m$ and all $i\leq n$, the $i^{\text{th}}$
coordinate $\eta^j_i$ of $\eta^j$ is either equal to $\chi_i$, or to
some large number number $N\gg 0$.  The Equation \eqref{polys}
still holds if we wiggle these large numbers a little bit,
hence the polynomial $P^{\eta^j}$ must be independent of the variable
$\eta^j_i$ whenever $\eta^j_i\neq \chi_i$.  Thus we may substitute
$\chi$ for each $\eta^j$ on the right-hand side, and we obtain the equation
\begin{equation*}
\PrA(\chi) = \sum_{j=1}^m\a_j P^{\eta^j}(\chi).
\end{equation*}
This equation clearly holds in a neighborhood of $\chi$,
hence we obtain an equation of polynomials
$$\PrA = \sum_{j=1}^m\a_j P^{\eta^j}.$$
This proves that $U\subs\R\big\{\PrA \mid A\text{ admissible}\big\}$, 
and thereby completes the proof of Theorem \ref{int}.
\end{proofint}

\begin{example}\em
Let us consider an example in which $n=4$ and $d=2$.  The picture on the left-hand side of 
Figure \ref{two} shows a polytope $\Dr$ where $\chi=(0,1,1,0)$, 
and the picture on the right shows $\Dr_{\{1,4\}}$.
\begin{figure}[h]
\begin{center}
\includegraphics[height=30mm]{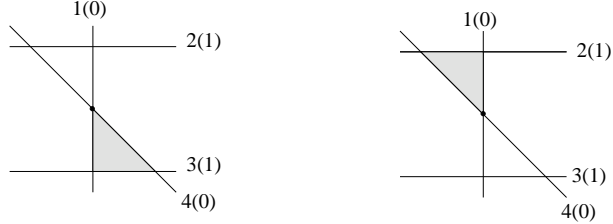}
\caption{The number outside of the parentheses denotes the index $i$ of the half space,
and the number inside denotes the value of $\chi_i$.  This value is equal (up to sign) to the distance from the
boundary of the $i^{\text{th}}$ half space to the origin of $\tddr$, which is marked with a black dot.}
\label{two}
\end{center}\end{figure}

\noindent
The key to the proof of Theorem \ref{int} is our ability
to express the characteristic function of $\Dr_{\{1,4\}}$ in terms of the characteristic
functions of $\D^{\eta^j}$ for some finite set $\{\eta^1,\ldots,\eta^m\}$ of simple integral vectors,
as we did in Equation \ref{express}.
Since $\{1,4\}$ has two elements, the procedure described
in Lemma \ref{key} must be iterated twice, and the result will have a total
of $2^2=4$ terms, as illustrated in Figure \ref{char}.
The first iteration exhibits $\mathbf{1}_{\Delta_{\{1,4\}}}$ as an element
of $W_{\{4\}}^{bd}$ by expressing it as the difference of the characteristic
functions of two (unbounded) regions.
With the second iteration, we attempt to express each of these two characteristic functions
as elements of $W_{\{1,4\}}^{bd} = W_{\{1,4\}}$.
This attempt must fail, because each of the two functions that we try
to express has unbounded support.  But the failures cancel out, and we succeed
in expressing the {\em difference} as an element of $W_{\{1,4\}}$.
\begin{figure}[!ht]
\begin{center}
\includegraphics[height=136mm]{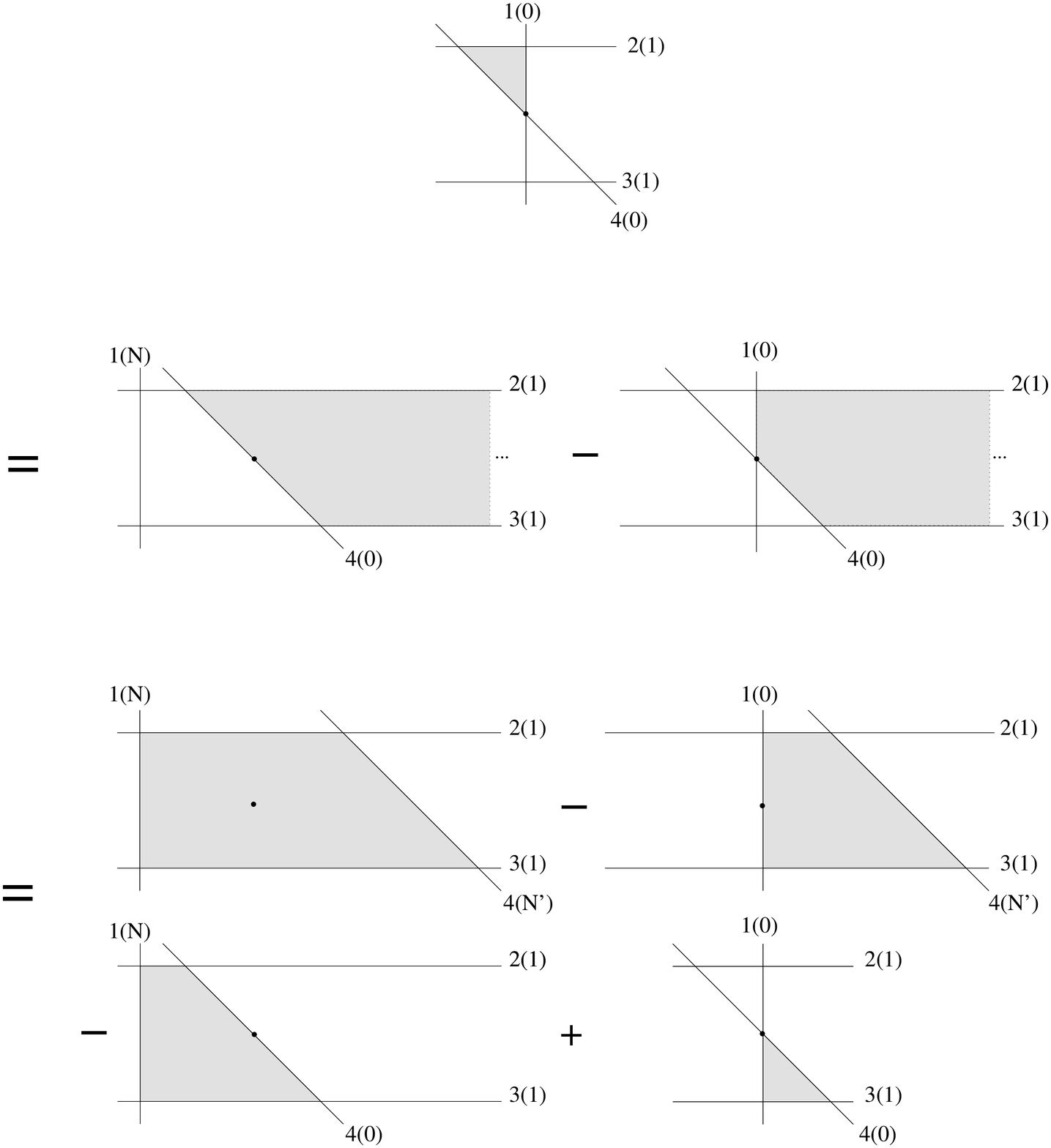}
\caption{An equation of characteristic functions.
The procedure of Lemma \ref{key} have produced two undetermined
large numbers, which we call $N$ and $N'$.}\label{char}
\end{center}\end{figure}

Let's see what happens when we take volume polynomials in the equation
of Figure \ref{char}.  The two polytopes on the top line have different
volumes, but the same volume polynomial, hence these two terms cancel.
We are left with the equation
$$P_{\{1,4\}}^{(0,1,1,0)} = P^{(0,1,1,0)} - P^{(N,1,1,0)},$$
which translates into
$$\half\(-\chi_1+\chi_2-\chi_4\)^2 = \half\(\chi_1+\chi_3+\chi_4\)^2 -
\(\chi_2+\chi_3\)\(\chi_1+\chi_4+\half \chi_3 - \half \chi_2\).$$
\end{example}
\end{section}

\paragraph{\bf Acknowledgments.}  This paper benefited greatly from conversations
with Matthias Beck and Rebecca Goldin, and from the work of William Crawley-Boevey.

\footnotesize{

}

\end{document}